\documentclass[11pt]{amsart}
\usepackage{amsmath}
\usepackage{amssymb}
\usepackage{epsfig}
\usepackage{a4wide}

\usepackage{pifont}

\newcommand{\B}{\ensuremath{\mathbb{B}}}

\newcommand{\Z}{\ensuremath{\mathbb{Z}}}
\newcommand{\R}{\ensuremath{\mathbb{R}}}
\newcommand{\T}{\ensuremath{{\mathcal T}}}
\newcommand{\F}{\ensuremath{{\mathcal F}}}
\newcommand{\bs}[1]{\boldsymbol{#1}}
\newcommand{\E}{\ensuremath{\mathbb{E}}}

\DeclareMathOperator{\inn}{int}
\DeclareMathOperator{\cl}{cl}

\DeclareMathOperator{\dens}{dens}
\DeclareMathOperator{\conv}{conv}

\newtheorem{thm}{Theorem}[section] 

\newtheorem{prop}[thm]{Proposition}
\newtheorem{defi}[thm]{Definition}

\begin{document}

\title{Duality of Model Sets generated by Substitutions}

\author{D.~Frettl\"oh} 

\address{Fakult\"at f\"ur Mathematik, Universit\"at
  Bielefeld, Postfach 100131, 33501 Bielefeld,  Germany}
\email{dirk.frettloeh@math.uni-bielefeld.de}
\urladdr{http://www.math.uni-bielefeld.de/baake/frettloe}

\maketitle

\begin{center}
{\em Dedicated to Tudor Zamfirescu \\
on the occasion of his sixtieth birthday}
\end{center}

\begin{abstract} 
The nature of this paper is twofold: On one hand, we will give a short
introduction and overview of the theory of model sets in connection
with nonperiodic substitution tilings and generalized Rauzy fractals.
On the other hand, we will construct certain Rauzy fractals and a
certain substitution tiling with interesting properties, and we will
use a new approach to prove rigorously that the latter one arises 
from a model set. 
The proof will use a duality principle which will be described in
detail for this example. This duality is mentioned as early as
1997 in \cite{gel} in the context of iterated function systems, but it 
seems to appear nowhere else in connection with model sets.
\end{abstract}

\section{Introduction\label{intro}}
One of the essential observations in the theory of nonperiodic, but
highly ordered structures is the fact that, in many cases, they can be
generated by a certain projection from a higher dimensional (periodic)
point lattice. This is true for the well-known Penrose tilings (cf.~
\cite{gs}) as well as for a lot of other substitution tilings, where
most known examples are living in $\E^1,\E^2$ or $\E^3$. The
appropriate framework for this arises from the theory of model sets
(\cite{mey}, \cite{moo}). Though the knowledge in this field has grown
considerably in the last decade, in many cases it is still hard to
prove rigorously that a given nonperiodic structure indeed is a model
set.  
 
This paper is organized as follows. Section \ref{sec:subst} briefly 
collects some well--known facts about substitution
tilings. These facts will be used later. Section
\ref{sec:modset} contains some basics about model sets, in particular
in connection with substitution tilings. In Section 1.3 we give the
definition of Rauzy fractals. This can be done completely in the
framework of model sets; it is just a question of the terminology.  
In Section \ref{sec:ansatz} we introduce a family of one--dimensional 
substitutions which are generalizations of the well--known Fibonacci
sequences. Four of these possess a substitution factor which is a
PV--number, so these are candidates for being model sets. Section
\ref{sec:n4} is dedicated to one of these substitutions. Therein, the
corresponding Rauzy fractal is obtained. Section \ref{sec:dual} shows
the construction of the dual substitution tiling, which is a
two--dimensional tiling with some interesting properties. Using the
duality, we will show that this dual tiling is a model set. Hence,
this is one of the few known cases of a two--dimensional model set, 
generated by a substitution, where the substitution factor is 
not an algebraic number of degree one or degree two (and the proof
being rigorous). Section \ref{sec:remarks} contains some additional
remarks, and further examples of dual substitutions, where the duality
principle can be used in order to prove that these are model sets, too.

Let us fix some notation. Throughout the paper, $\cl(M)$ denotes the
closure of a set $M$, $\inn(M)$ the interior of $M$, $\partial M$ the
boundary of $M$, and $\#M$ the cardinality of $M$. $\E^d$ denotes the
$d$--dimensional Euclidean space, i.e., $\R^d$ equipped with the
Euclidean metric $\| \cdot \|$. $\B^d$ denotes the unit ball $\{ {\bs
  x} \, | \, \|{\bs   x}\|\le 1 \}$. If not stated otherwise, $\mu$
denotes the $d$--dimensional Lebesgue measure, where $d$ will be clear
from the context.    
A \emph{tile}  is a nonempty compact set $T \subset \E^d$ with the property
$\cl(\inn(T))=T$. A \emph{tiling} of $\E^d$ is a collection of tiles 
$\{T^{}_n\}^{}_{n \ge 0}$ which covers $\E^d$ and contains no 
{\em overlapping} tiles, i.e., $\inn(T^{}_k) \cap \inn(T^{}_n) =
\varnothing$ for $k \ne n$. In other words: A tiling is both a
covering and a packing of $\E^d$. A tiling $\T$ is called
\emph{nonperiodic}, if the only solution of  $\T+\bs{x}=\T$  is
$\bs{x}=\bs{0}$. 

\subsection{Substitution Tilings} \label{sec:subst}
\begin{figure}[h] 
\epsfig{file=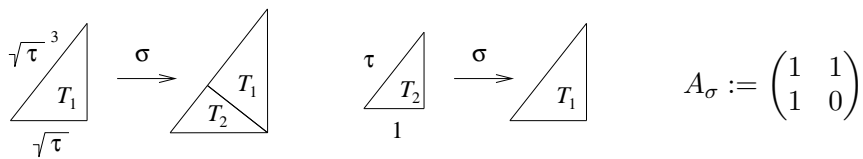} \hspace{10mm}
\begin{minipage}[b]{30mm}
$A_{\sigma}:=\begin{pmatrix} 1 & 1 \\ 1 &  0 \end{pmatrix}$
\vspace*{5mm}
\end{minipage}
\caption{'Golden triangles'. A substitution rule in $\E^2$ with
  two prototiles $T_1,T_2$ and substitution factor $\lambda=\sqrt{\tau}$,
  where $\tau=\frac{\sqrt{5}+1}{2}$ is the golden mean. 
 \label{inflbsp}}
\end{figure}

A convenient method to generate nonperiodic tilings is by a 
{\em substitution}: Choose a set of {\em prototiles}, that is, a set of
tiles $\F:=\{T_1, T_2, \ldots, T_m\}$. Choose a 
{\em substitution factor} $\lambda \in \R, \lambda > 1$ and a rule how to
dissect $\lambda T_i$ (for every $1 \le i \le m$) into tiles, such that
any of these tiles is congruent to some tile $T_j \in \F$. It will be
sufficient for this paper to think of a substitution like in
Fig.~\ref{inflbsp}. For completeness, we give the proper definition
here.

\begin{defi} \label{substdef}
Let $\F:=\{T_1, T_2, \ldots, T_m\}$ be a set of tiles, the 
{\em  prototiles}, and $\lambda>1$ a real number, the 
{\em substitution factor}.  Let $\lambda T_i =  \varphi_1(T_{i_1}) \cup
  \varphi_2(T_{i_2}) \cup \ldots  \cup \varphi_{n(i)}(T_{i_{n(i)}})$,  
such that every $\varphi_j$ is an isometry of $\E^d$ and the
involved tiles do not overlap. Furthermore, let 
$\sigma(\{T_i\})= \{ \varphi_1(T_{i_1}), \varphi_2(T_{i_2}), \ldots,
\varphi_{n(i)}(T_{i_{n(i)}}) \}$. Let ${\bs S}$ be the set of all
configurations of tiles congruent to tiles in $\F$. That is, ${\bs S}$ is
the set of all sets of the form  $\{T_{j_i} +{\bs x}_i \, | \, T_{j_i} \in
\F, i \in I, {\bs x}_i \in \E^d \}$, where $I$ is some 
finite or infinite index set. By the requirement
$\sigma(\{T_i+{\bs x}\})=\sigma(\{T_i\})+\lambda {\bs x}$ (${\bs x}
\in \E^d$), $\sigma$ extends in a 
unique way to a well--defined map from ${\bs S}$ to ${\bs S}$. 
Then $\sigma: {\bs S} \to {\bs S}$ is called a {\em substitution}. 
Every tiling of $\E^d$, where any finite part of it is congruent to a
subset of some $\sigma^k(T_i)$ is called a 
{\em substitution tiling (with  substitution $\sigma$)}.
\end{defi}  

More general, this definition makes sense if we replace $\lambda$
by some expanding linear map. That is, replace $\lambda$ in the
definition above by a matrix $Q$ such that all eigenvalues of $Q$ are
larger than 1 in modulus. In this full generality, many of the following
results may break down. Anyway, in Section \ref{sec:n4} we will use the
generalization in one certain case, where it will cause no problem. 

To a substitution $\sigma$ we assign a {\em substitution matrix}:
$A_{\sigma}:=(a_{ij})_{1 \le i,j \le m}$, where $a_{ij}$ is the number
of tiles $T_j+{\bs x}$ (tiles of type $T_j \in \F$) in $\sigma(\{T_i\})$. 
In this paper we will consider substitutions with {\em primitive}
substitution matrices only. Recall: A nonnegative matrix $A$ is called 
{\em primitive}, if some power $A^k$ is a strictly positive matrix. In
this situation we can apply the following well--known theorem.

\begin{thm}[Perron--Frobenius] \label{perron}
Let $A \in \R^{m \times m}$ be a primitive nonnegative matrix. 
Then $A$ has an eigenvalue $0< \lambda \in \R$, which is simple and
larger in modulus than every other eigenvalue of $A$. This
eigenvalue is called {\em Perron--Frobenius--eigenvalue} or shortly 
{\em  PF--eigenvalue}.    \\
The corresponding eigenvector can be chosen such that it is strictly
positive. Such an eigenvector is called {\em PF--eigenvector} of $A$.   
No other eigenvalue of $A$ has a strictly positive eigenvector.
\end{thm} 

Using this theorem, it is a simple exercise to show the following
facts (cf.~\cite{fogg}, Thm.~1.2.7): We consider the 
{\em relative frequency} of the tiles of type $T_i$ in a tiling $\T$,
i.e., the ratio of the number of tiles of type $T_i$ in $\T$ and
the number of all tiles in $\T$. Formally, this is  defined by  
$\lim_{r\to \infty} \frac{\# \{ T_i + {\bs x} \, \in \T \, | \, {\bs x} 
  \in r\B^d  \}}{\# \{ T+{\bs x} \, \in \T \, | \, {\bs x} \in r\B^d,
  \, T   \in \F \}}$. Since we
consider primitive substitution tilings only, it is true that this
limit exists. Moreover, this limit exists {\em uniformly}, i.e., for
any translate  $\T+{\bs x}$ we will obtain the same limit.
\begin{prop} \label{freqvol}
Let $\sigma$ be a primitive substitution in $\E^d$ with substitution
factor $\lambda$ and prototiles $T_1, \ldots, T_m$. Then the
PF--eigenvalue of a $A_{\sigma}$ is 
$\lambda^d$. In particular, a substitution factor is always an algebraic
number. \\ The normalized (right) PF--eigenvector ${\bs v}= 
(v_1,\ldots,v_m)^T$ of $A_{\sigma}$ contains the relative frequencies
of the tiles of  different types in the tiling in the following sense:
The entry $v_i$ is the relative frequency of $T_i$ in $\T$.\\ 
The left PF--eigenvector
(resp.~the PF--eigenvector of $A^T$) contains the $d$--dimensional
volumes of the different prototiles, up to scaling.
\end{prop}

Some people don't like the term 'left eigenvector'. Those may
replace 'left eigenvector of $A$' by 'eigenvector of $A^T$, written as
a row vector' wherever it occurs. 
In the next section, we will associate with a tiling a related point
set: By replacing every prototile $T_i$ with a point in its
interior, a tiling $\T$ gives rise to a {\em Delone set} $V_{\T}$. 
(Note that this defines $V_{\T}$ not uniquely. It depends on the
distinct choice of points in the prototiles. But the following is
true for any Delone set $V_{\T}$ constructed out of a
primitive substitution tiling $\T$ in the described way.) 
A Delone set is a point set $V \subset \E^d$ which is  
{\em uniformly discrete} (there exists $r>0$ such that $\# (V \cap
r\B^d + {\bs x}) \le 1$ for all ${\bs x} \in \E^d$)
and {\em relatively dense} (there exists $R>0$ such
that $\# (V \cap R\B^d + {\bs x}) \ge 1$ for all ${\bs x} \in \E^d$).    
Since the frequencies of the prototiles tiles (cf. remark after Theorem
\ref{perron}) are well--defined for all tilings $\T$ considered in this
paper, it makes sense to define the density of the associated
Delone set $V_{\T}$: $\dens(V_{\T})=\lim_{r\to \infty} \frac{\# V_{\T}
  \cap r\B^d} {\mu( r \B^d)}$. It is not hard to show the following
identity. 

\begin{prop} \label{densfreqvol}
Let $\T$ be a substitution tiling with substitution
$\sigma$, let ${\bs v}$ be the right normalized PF--eigenvector of
$A_{\sigma}$ (containing the relative frequencies of the prototiles),
let ${\bs w}^T$ be the left PF--eigenvector of $A_{\sigma}$, such that 
${\bs w}^T$ contains the $d$--dimensional volumes of the prototiles, and
$V_{\T}$ a Delone set constructed out of $\T$ by replacing every
prototile with a point in its interior. Then
\[ \dens(V_{\T}) = \left( {\bs w}^T \cdot{\bs v} \right)^{-1} \]     
\end{prop}

{\em Example:} The golden triangle substitution. The substitution
rule, together with substitution factor, substitution matrix etc. is
shown below. 

This substitution uses two prototiles $T_1,T_2$. 
The PF--eigenvalue of the substitution matrix $A_{\sigma}$ 
is  $\tau=\lambda^2$. The normalized PF--eigenvector is
${\bs v}:=(\tau^{-1}, \tau^{-2})^T$. Thus, the relative
frequency of the tiles of type $T_1$ (resp.~$T_2$) in the tiling is
$\tau^{-1}$ (resp.~$\tau^{-2}$). $A_{\sigma}$ is symmetric, thus 
${\bs v}^T$ is a left PF--eigenvector. The areas (2--dim volumes) of the
two prototiles are $\mu(T_1)=\sqrt{\tau}^3, \mu(T_2)=\sqrt{\tau}$. If
we replace every triangle in a corresponding substitution tiling by
a point in its interior, then by Prop.~\ref{densfreqvol} the density
of the point set is $\left( (\sqrt{\tau}^3,\sqrt{\tau})\cdot(
  \tau^{-1}, \tau^{-2})^T \right)^{-1} =
\frac{\tau^{3/2}}{(\tau^2+1)}=0.56886448\ldots$

\subsection{Model Sets} \label{sec:modset}
Model sets are a special kind of Delone sets. They show a high degree
of local and global order. For example, a model set $V$ and its
difference set $V-V := \{ {\bs x} - {\bs y} \, | \, x,y \in V \}$
do not differ too much in the following sense (cf.~\cite{mey}): There
is a finite set $F$ such that 
\[ V-V \subseteq V+F. \]
Lattices do have this property with $F=\{0\}$, and indeed lattices are
special cases of model sets.  Unfortunately, there is no reference
known to the author which gives a comprehensive overview about model
sets. A good reference is \cite{moo}. A comprehensive
overview for the case of one--dimensional tilings is contained in
\cite{fogg}, mainly chapter 7 and 8. Therein, many connections to
combinatorics, ergodic theory and number theory are compiled (but the
term 'model set' cannot be found in this book, the term 'geometric
representation' is used instead). 

A model set is defined by a collection of spaces and maps (a
so--called cut--and--project scheme) as follows.
In (\ref{cps}), let $\Lambda$ be a lattice of full rank in
$\E^{d+e}$, $\pi_1,\pi_2$ projections such that  
$\pi_1 |_{\Lambda}$ is injective, and $\pi_2(\Lambda)$ is dense in
$\E^e$. Let $W$ be a nonempty compact set --- the so--called 
{\em window  set} --- with the properties $\cl(\inn(W))=W$
and $\mu(\partial W)=0$, where $\mu$ is the Lebesgue measure in $\E^e$.   
\begin{equation} \label{cps}
\begin{array}{ccc} 
\E^d & \stackrel{\pi_1}{\longleftarrow} \E^{d+e}
  \stackrel{\pi_2}{\longrightarrow} & \E^e \\
\cup & \cup & \cup \\
V & \Lambda & W
\end{array} 
\end{equation}
Then the set $V= \{ \pi_1({\bs x}) \, | \, {\bs x} \in \Lambda,
\pi_2({\bs x}) \in W \}$  
is called a {\em model set}\footnote{To be precise, this defines a
  {\em regular} model set. Without the property $\mu(\partial W)=0$,
  $V$ is still a model set, but not a regular one.}.
Obviously, $V$ is a set of points in
$\E^d$. A tiling can be obtained by applying certain rules to $V$,
e.g., considering the Voronoi cells of $V$. In the case $d=1$ it is
quite simple to get a tiling: If $V=\{ {\bs x}_i \, | \, i \in \Z \}$,
such that ${\bs x}_i<{\bs x}_{i+1}$, then $\T= \{ [{\bs x}_i, {\bs
  x}_{i+1} ] \, | \, i \in \Z \}$ is a tiling of $\E^1$. 

Since $V \subset \pi_1(\Lambda)$ and $\pi_1|_{\Lambda}$ is injective, 
every ${\bs x} \in V$ corresponds to exactly one point $\pi_1^{-1}({\bs x}) \in
\Lambda$, and thus exactly to one point $\pi_2(\pi_1^{-1}({\bs x})) \in W$.
Therefore, we are able to define the {\em star map}:
\begin{equation} \label{starmap}
  \ast : \pi_1(\Lambda) \to \E^e, \qquad {\bs x}^{\ast}  =
\pi_2(\pi_1^{-1}({\bs x})) 
\end{equation}
This allows for a convenient notation. E.g., we have
\begin{equation} \label{clvstern}
\cl(V^{\ast})=W, 
\end{equation}
where $V^{\ast} = \{ {\bs x}^{\ast} \, | \, {\bs x} \in V \}$. 
 
{\em Remark:} In general, $\E^d$ and $\E^e$ in (\ref{cps}) can be
replaced by some locally compact abelian groups $G$ and $H$. Then
$\mu$ denotes the Haar measure of $H$. This is necessary to show that,
for example, the chair tiling is a model set (cf~\cite{bms}). 
If the substitution matrix $A_{\sigma}$ of a given substitution tiling
is {\em unimodular} (i.e., $\det(A_{\sigma} = \pm 1$), then it
suffices to consider Euclidean spaces only. If not, then one needs to
consider more general groups (cf.~\cite{sing}). In this paper, all
occurring substitution matrices are unimodular, therefore we will 
consider Euclidean spaces only.   

Given a substitution tiling, one may ask if this tiling can be
obtained as a model set $V$. By the following theorem, this is only
possible if the substitution factor is a PV number or a Salem number.    

\begin{defi}
Let $\lambda \in \R$ be an algebraic integer. $\lambda$ is called a 
{\em  Pisot Vijayaraghavan number} or shortly {\em PV number}, 
if $|\lambda|>1$, and for all its algebraic conjugates $\lambda_i$
holds: $|\lambda_i|<1$. \\
If for all conjugates $\lambda_i$ of $\lambda$ holds
$|\lambda_i|\le 1$, with equality for at least one $i$,  then
$\lambda$ is called a {\em Salem number}.  
\end{defi}

\begin{thm}[Meyer] \label{pvsatz}
If $V$ is a model set, $1 < \lambda \in \R$ and if $\lambda V \subseteq V$,
then $\lambda$ is a PV number or a Salem number. 
\end{thm}

It is known that every substitution $\sigma$ gives rise to at least
one tiling $\T$ with the property $\sigma(\T)=\T$. So,  
if $V$ is a point set derived from $\T$ (where the substitution factor 
is $\lambda$), say, $V$ being the vertex set of this tiling, then by the
definition of a substitution follows $\lambda V \subseteq V$.  
Now, if $V$ is also a model set, by the theorem above the substitution
factor must be a PV number or a Salem number.

The remaining things to check in order to show that $V$ is a model set
are the following (cf.~(\ref{cps})):
\begin{itemize}
\item Determine the lattice $\Lambda$,
\item determine the window set $W$,
\item show $\mu(W)>0$,
\item show $\mu(\partial W)=0$,
\item show $\dens(V) = \dens(\{ \pi_1({\bs x}) \, |
  \, {\bs x} \in \Lambda, \pi_2({\bs x}) \in W \})$. 
\end{itemize} 
There are standard constructions for obtaining
$\Lambda$ and $W$. But the latter three items are hard to decide
rigorously in general, or even in special cases. For an example how
hard it may be to give a proper proof, cf.~\cite{bs}. In the next section,
we will use a duality principle to prove that the three latter items
are fulfilled for the considered tiling, thus proving that it arises 
from a model set (i.e., the associated point set is a model set). But
first, we will discuss the above list in more detail. 

Starting with a substitution tiling, or better a substitution point
set $V$ in $\E^1$ (cf.~the remark after (\ref{cps})), the standard
construction for the lattice is as follows: Let $\lambda$ be the
substitution factor (where $\lambda$ from now on is always assumed to be
a PV number) and $\lambda_2,\ldots,\lambda_m$ its algebraic conjugates. 
Then \[ \Lambda= \langle (1,1,\cdots,1)^T,
(\lambda,\lambda_2,\cdots,\lambda_m)^T,
(\lambda^2,{\lambda_2}^2,\ldots,{\lambda_m}^2)^T, \ldots,   
(\lambda^{m-1},{\lambda_2}^{m-1},\ldots,{\lambda_m}^{m-1})^T \rangle_{\Z}. \]
Now, $\pi_1$ (resp.~$\pi_2$) can be chosen as the canonical projections
from $\E^m \to \E^1$ (resp.~$\E^m \to \E^{m-1}$). To be precise, if 
${\bs  x}=(x_1,x_2,\ldots,x_m)^T$, then $\pi_1({\bs x})=x_1,
\pi_2({\bs   x})=(x_2,\ldots, x_m)^T$.   
At this point, we are already able to compute the window set
$W$ numerically: Since the left PF--eigenvector 
of the substitution matrix contains the the lengths of the 
prototiles (cf.~remark after Thm.~\ref{perron} above), the
lengths can be chosen such that $V \subset \Z[\lambda]$. That is, every
${\bs x} \in V$ can be expressed as ${\bs x}=\alpha_0+\alpha_1 \lambda+
\cdots+ \alpha_{m-1} \lambda^{m-1}$, where $\alpha_i \in \Z$.   
Then 
\[ {\bs x}^{\ast} = (\alpha_0+\alpha_1 \lambda_2+ \cdots+ \alpha_{m-1}
{\lambda_2}^{m-1}, \cdots, \alpha_0+\alpha_1 \lambda_{m-1}+ \cdots+
\alpha_{m-1} {\lambda_{m-1}}^{m-1})^T. \]
More general, we may apply the star map to $V$, obtaining the
(countable) set $V^{\ast}$, which in turn yields the window set
$W=\cl(V^{\ast})$. Actually, some of the pictures in the next section
were created in essentially in this way: Using the substitution to
generate a large finite subset $F \subset V$, where the points are
expressed as ${\bs x}$ above, we plot all points of $F^{\ast}$. 
If this set $F^{\ast}$ appears to be bounded, even if the number of
points is increased, then this is an indication that the considered
set $V$ is a model set. 

To prove rigorously that $V$ is actually a model set, one proceeds by
showing $\mu(W)>0$ and $\mu(\partial W)=0$. If $W$ is a polyhedron,
then this is an easy task. But in the generic case $W$ is of fractal
nature, and it may be even hard to prove that $W$ is of positive
measure. 

The last thing to prove is $V=\{ \pi_1({\bs x}) \, | \, {\bs x} 
\in \Lambda, \pi_2({\bs x}) \in W \}$, and not just a too small
subset. More precisely: Let $U:=\{ \pi_1({\bs x}) \, | \, {\bs x}
\in\Lambda,\pi_2({\bs x})\in W\}$. Then it must hold 
$\dens(V)=\dens(U)$. In this context the following
theorem turns out to be helpful \cite{sch}.
  
\begin{thm} \label{densvol}
Let $V$ be a model set given by (\ref{cps}). Then
\[ \dens(V) = \mu(W) / \det \Lambda, \]
where $\det \Lambda$ denotes the determinant of the generator matrix of
$\Lambda$. Equivalently, $\det \Lambda$ is the volume of any measurable
fundamental domain of $\Lambda$.   
\end{thm}

\subsection{Rauzy Fractals} \label{sec:rauzy}
As mentioned above, 'Rauzy fractal' 
is essentially just another term for the window set $W$. The
term 'Rauzy fractal' occurs in the context of dynamical systems and
number theory, while the terms 'window' or 'window set' arose from the
theory of quasicrystals. Originally, 'Rauzy fractal' denoted one special
fractal of this kind, namely the one arising from the substitution 
\[ \sigma(1)=1\;2, \quad  \sigma(2) =   1\; 3, \quad \sigma(3)= 1, \]
which is analyzed in \cite{r}. Therefore, the term 'generalized Rauzy
fractal' is sometimes used to 
emphasize the fact that one deals with the more general case. 
(Of course, this is just a symbolic substitution, but it is trivial to
formulate the geometric substitution according to Def. \ref{substdef},
cf.~comment after (\ref{substn})).  

According to \cite{fogg}, the definition of a Rauzy fractal is
essentially as follows: if $\sigma$ is a substitution for $m$
prototiles in $\E^1$, where the factor is a PV number of algebraic
degree $m$, then the Rauzy fractal is just the window set $W$ as in
(\ref{clvstern}), where one does not care about $\mu(W)>0$ or
$\mu(\partial W)=0$ or Theorem \ref{densvol}. Note that therefore it
is much more difficult to show that a given substitution point set is
a model set than it is to compute the corresponding Rauzy fractal.       

Before we proceed, let us mention a theorem which is important in this
context and which is necessary in the next section. 

\begin{thm}[\cite{hut}] \label{hutch}
Let $(X,d)$ be a complete metric space and $\{f_1,\ldots,f_n\}$ a
finite set of contractive maps (i.e., $\exists c<1: \forall x,y \in X:
d(f(x),f(y)) \le c\, d(x,y)$). Then there is a unique compact set $K
\subset X$ such that $K = \bigcup_{i=1}^n  f_i(K)$. 
\end{thm}

This theorem is a consequence of the Banach Fixed Point Theorem, and is
easily generalized to partitions of $K$: Let $\{f_1,\ldots,f_n\}$ be
as above and $I_j \subseteq \{1,\ldots,n\}$ for $1 \le j \le \ell$,
then there is a unique tuple $(K_1,\ldots,K_{\ell})$ of compact sets
$K_j$, such that $K_j = \bigcup_{i \in I_j}  f_i(K_{j_i})$ for all 
$1 \le j \le \ell$ (\cite{kli}). In the next section, this theorem will
help to show the compactness of the window set $W$.

\section{Two Rauzy Fractals and the Duality Principle} \label{sec:zwo}
We will describe the duality principle by choosing a certain
one--dimensional substitution, constructing the corresponding 
Rauzy fractal and the associated dual tiling, and then we will
use the duality principle to show that the dual tiling is a model set.  

\subsection{Some Rauzy Fractals} \label{sec:ansatz}
It is easy to see, that for any nonnegative integer $n\times n$--matrix
$A$ there is a substitution tiling in $\E^1$ with substitution matrix
$A$. Actually, from any such matrix one can immediately deduce a
substitution, since there are only few geometric restrictions for
tilings in one dimension. (The same problem for higher
dimensions becomes quite difficult!) Here we consider the following
matrices for $n \ge 2$: 
\[ M_n:=(m_{ij})_{1 \le i,j \le n}, \quad \mbox{where} \; m_{ij} = \left\{
    \begin{array}{lcl} 1 & : & i+j >n \\ 0 & : & i+j \le n 
\end{array} \right. \]
Obviously all these matrices are primitive.
A substitution with substitution matrix $M_n$ for the prototiles
$1,2,\ldots,n$ is given by 
\begin{equation} \label{substn}
\begin{array}{ccr} 1 & \to & n \\ 2 & \to & n\!-\!1 \; n \\ & \vdots & 
\vdots \\ n & \to & 1 \; 2 \cdots n\!-\!1 \; n \end{array} 
\end{equation}   
Note that this is just a symbolic substitution. But in one dimension
the geometric realization is obtained from this scheme in a unique way:
The prototiles are intervals of different lengths, and the lengths are
given by the PF--eigenvector of $M_n$. There is a nice trick to deduce the
substitution factor directly from this substitutions: For fixed 
$n \in \Z, n \ge 2$, let $s_k:=\sin(\frac{k \pi}{2n+1})$.  
Because of the wonderful formula (see \cite{nida})
\begin{equation} \label{wformel}
 \frac{s_k}{s_1} s_i =
\sum_{\nu=0}^{k-1}s_{i+1-k+2\nu} 
\end{equation}  
it follows $\frac{s_n}{s_1} s_1 = s_n,  \frac{s_n}{s_1} s_2 = s_{n-1}+s_n,
\ldots, \frac{s_n}{s_1} s_n = s_1 + s_2 + \cdots + s_n$. Therefore the
substitution factor (which is the PF--eigenvalue of $M_n$) is
$\frac{s_n}{s_1}$, and the length of tile $i$ is $s_i$. 
Now, we determine which of these substitutions 
define a Rauzy--fractal. The first condition to be fulfilled is that
the substitution factor must be a PV--number\footnote{Or a Salem
  number, but these don't play a role here, so we don't mention them
  explicitly. }, cf.~Theorem 
\ref{pvsatz}. This is the case for $n=2,3,4$ and $7$, and for no
other $n \le 30$. (This has been checked by a Maple program, and we
don't expect any more PV--numbers as PF--eigenvalues of $M_n$ for
larger values of $n$.) The corresponding characteristic polynomials of
$M_n$ are:
\[ x^2-x-1, \; x^3-2x^2-x+1, \; x^4-2x^3-3x^2+x+1, \;
x^7-4x^6-6x^5+10x^4+5x^3-6x^2-x+1. \] 
The latter two are reducible over \Z. The irreducible polynomials
defining $\frac{s_2}{s_1}\, (n=2),\, \frac{s_3}{s_1}\, (n=3), \,
\frac{s_4}{s_1} \, (n=4)$ and $\frac{s_7}{s_1} \, (n=7)$ are 
\begin{equation}\label{charpoly}
 x^2-x-1, \; x^3-2x^2-x+1, \; x^3-3x^2+1, \; x^4-4x^3-4x^2+x+1. 
\end{equation}
Therefore, the algebraic degree of these numbers is 2,3,3,4, resp. 
Here we will focus mainly on the third case ($n=4$). The first case
($n=2$) leads to the well known Fibonacci sequences (see Section
\ref{sec:remarks}, or \cite{fogg}, Sections 2.6 and 5.4). 
The fourth case  
($n=7$) will lead to a  three-dimensional Rauzy fractal in the end. 
The third case ($n=4$) will show the methods and the typical problems
which may occur, and therefore we discuss this case in the following
in detail. This part may be regarded as a recipe how to compute other
Rauzy fractals explicitly. 

\subsection{The case ${\bs n=4}$} \label{sec:n4}
In this case we have got four prototiles. Since we consider one--dimensional
tilings here, the tiles can be chosen to be intervals. We denote the
four prototiles by $L$ ('large'), $M$ ('medium'), $S$ ('small') and
$X$ ('extra small'). Their lengths are already known: They are given
by a left PV--eigenvector of the substitution matrix, or by
(\ref{wformel}), so the values can be chosen as 
$\frac{s_4}{s_1}, \frac{s_3}{s_1}, \frac{s_2}{s_1},
\frac{s_1}{s_1}=1$, resp.~  (or any positive multiples of these
values).  In order to construct the Rauzy fractal, we need this values 
as expressions in $\lambda=\frac{s_4}{s_1}$ (cf. Section \ref{sec:modset}). 
W.l.o.g, we choose $X:=[0,1]$. The substitution now reads
\begin{equation}\label{subst}
 X \to L, \; S \to ML, \; M \to SML, \; L \to XSML. 
\end{equation}
Since $\lambda X = L$ it follows $L=[0,\lambda]$. Let $s,m,\ell$ denote
the lengths of $S,M,L$, resp. Since $\ell=\lambda$, we conclude from
(\ref{subst}) the equations
\[ \begin{array}{rcrcl} 
\lambda s & = &  m+ \lambda & & \\ 
\lambda m & = & s+m+ \lambda  & = & s(\lambda+1) \\ 
\lambda^2 & = & 1+s+m+\lambda & = & s(\lambda+1) +1 = \lambda m +1 
\end{array} \]   
and therefore $m=\frac{\lambda^2-1}{\lambda}$. From (\ref{charpoly})
follows $\lambda^3=3\lambda^2-1$, hence $\lambda^{-1}=3 \lambda -
\lambda^2$, and we obtain
\[ m = (\lambda^2-1)\cdot(3\lambda-\lambda^2) = \lambda^2 - 2 \lambda. \] 
Consequently,
\[ s = \lambda^2-\lambda-m-1 = \lambda - 1. \]
To apply the standard construction of Rauzy fractals, the number of
prototiles must equal the algebraic degree of the substitution
factor. So we need to get rid of one letter. By (\ref{subst}), a
tile $X$ in the tiling is always followed by a tile $S$. Since $x+s=\lambda$,
we can merge every pair $XS$ into one tile $L$. It is easy to check
that the new substitution 
\begin{equation} \label{subst2}
S \to ML, \; M \to SML, \; L \to LML 
\end{equation}
gives rise to the same tilings as (\ref{subst}) does, up to merging
pairs $XS$ into $L$. We identify these tilings with a point set
$V$, where $V$ contains all (right) endpoints of the tiles. If $0 \in
V$, then $V \subset \Z[\lambda]$, and every point $x \in V$ can be
expressed by $x=a + b \lambda + c \lambda^2$ ($a,b,c \in \Z$).  
The lattice $\Lambda \in \E^3$ now is given by 
\begin{equation} \label{mlambda}
\Lambda = \langle (1,1,1)^T, (\lambda, \lambda_2, \lambda_3)^T, (\lambda^2,
\lambda_2^2, \lambda_3^2)^T \rangle_{\Z}, 
\end{equation}  
where $\lambda_2, \lambda_3$ are the algebraic conjugates of $\lambda$. 
From $\lambda \lambda_2 \lambda_3 = -1$ and the fact that the values $s_k$
occur frequently above, one easily deduces
$\lambda_2=-\frac{s_1}{s_2}, \lambda_3=\frac{s_2}{s_4}$. 
This allows us to write down the star map now (cf.~(\ref{starmap})). 
Let $0 \in V$. If $x \in V$, then $x= a + b \lambda + c \lambda^2$ for some
$a,b,c \in \Z$, and 
\begin{equation} \label{ast}
x^{\ast} = a \begin{pmatrix} 1 \\ 1 \end{pmatrix} 
+ b \begin{pmatrix} \lambda_2 \\ \lambda_3 \end{pmatrix} 
+ c \begin{pmatrix} {\lambda_2}^2 \\ {\lambda_3}^2 \end{pmatrix}. 
\end{equation}
By (\ref{clvstern}), the window set $W$ is just the closure of $V^{\ast}$. 
A short computation (using $\lambda_i^3= 3 \lambda_i^2-1$ for $i=2,3$)
shows 
\[ (\lambda x)^{\ast} = \begin{pmatrix} \lambda_2 & 0 \\ 0 & \lambda_3
\end{pmatrix} x^{\ast}  =: Qx^{\ast}. \]

In order to show the compactness of $W$, and possibly to draw further
conclusions about $W$ (which is the Rauzy fractal wanted), we consider
a substitution set $V$ such that 
$\sigma(V)=V$. (As mentioned above, the considered substitution $\sigma$
can be easily formulated for the point set $V$ as well as for the
tiling $\T$.) For example, let
\begin{equation} \label{sv=v}
V = \cdots SMLLM\underline{L}LMLSMLLMLMLSMLLMLLMLSMLL \cdots 
\end{equation}
where the right endpoint of the underlined $L$ is located at zero.
(Note, that the condition $\sigma(V)=V$ together with (\ref{sv=v})
defines $V$ uniquely.) Denote by $V_L$ ($V_M,V_S$, resp.) the set of
all points in $V$ which are right endpoints of intervals $L$ ($M,S$,
resp.). Then obviously 
\[ V=V_L \cup V_M \cup V_S. \]
Moreover, from (\ref{subst2}) and $\sigma(V)=V$ follows
\begin{equation} \label{vgleich}
\begin{array}{lcl}
V_L & = & \lambda V \quad \cup \quad \lambda V_L -\lambda^2+\lambda, \\
V_M & = & \lambda V -\lambda, \\
V_S & = & \lambda V_M - \lambda^2+\lambda.
\end{array}
\end{equation} 
This can be seen as follows: One tile of type $M$ and one tile
of type $L$ occur in {\em every} substituted tile $T$, not depending
on the type of $T$. If $T=[a,b]$, then $\sigma(T)$ contains tiles
$[\lambda b - \lambda - m, \lambda b - \lambda]$ (of type $M$) and 
$[\lambda b -\lambda, \lambda b ]$ (of type $L$). Therefore, from 
$b \in V$ follows $\lambda b - \lambda \in V_M$,
yielding the second equation, and $\lambda b \in V_L$, yielding the first
part of the first equation. The other parts of (\ref{vgleich}) are
obtained analogously.

Applying the star map to (\ref{vgleich}), we obtain the following
equations for $V^{\ast}$:
\begin{equation} \label{vstern}
\begin{array}{lcl}
V^{\ast} & = & V_L^{\ast} \quad \cup \quad V_M^{\ast} \quad \cup 
\quad V_S^{\ast} \\
& & \\
V_L^{\ast} & = & Q V^{\ast} \quad \cup \quad Q V_L^{\ast}
-  ({\lambda_2}^2,  {\lambda_3}^2)^T + ( \lambda_2 , \lambda_3 )^T \\
V_M^{\ast} & = & Q V^{\ast} - ( \lambda_2 , \lambda_3)^T  \\
V_S^{\ast} & = & Q V_M^{\ast} - ( {\lambda_2}^2 , {\lambda_3}^2 )^T + 
(\lambda_2 ,   \lambda_3 )^T  
\end{array} 
\end{equation}
Obviously, $Q$ defines a contracting linear map, thus all maps
appearing in (\ref{vstern}) are contracting. Therefore, Theorem
\ref{hutch} applies: There is one unique compact $W=W_L \cup W_M \cup 
W_S$ which fulfills (\ref{vstern}). This assures the existence of this
set, so we obtained a Rauzy fractal $W$. 

\begin{figure}[h]
\epsfig{file=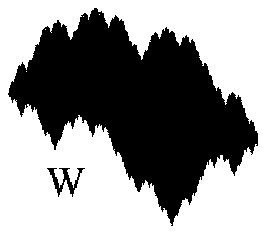} 
\epsfig{file=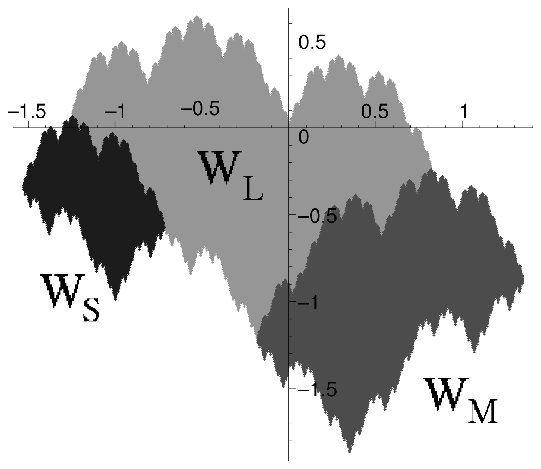} \hspace{5mm}
\epsfig{file=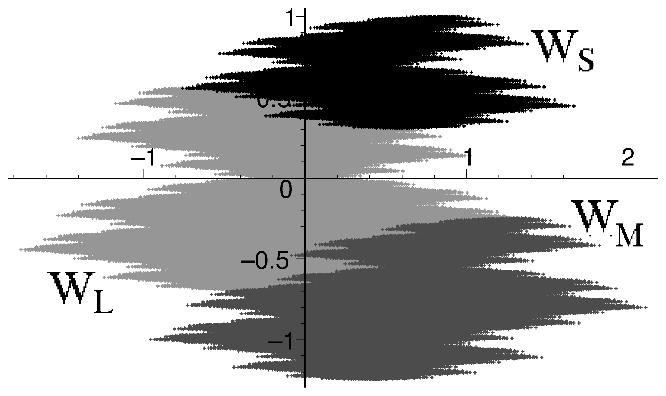}
\caption{Left: The Rauzy fractal $W$ arising from the substitution
  (\ref{subst2}). Middle: the same set, partitioned according to the
  distinction of tiles $S$, $M$ and $L$. Right: The Rauzy
  fractal arising from the substitution (\ref{substn}) for $n=3$:  $S
  \to L, \; M \to ML, \, L  \to SML$. \label{wkp9} } 
\end{figure}

To visualize $W$ and its fractal
nature, $W$ can be approximated numerically. In fact, we just produce a  
finite subset of $V$ and apply the star map to it. The resulting
points are plotted into a diagram, which is shown in Fig.~\ref{wkp9}
(left). If we distinguish the points in the diagram according to their
type ($S,M$ or $L$), the result is shown in Fig.~\ref{wkp9} (middle). 
The set $W$ is the union of the three sets $W_L, W_M, W_S$, where
$W_i=\cl(V_i^{\ast})$. From (\ref{vstern}) follows that $W_S$ is
an affine image of $W_M$, but not similar to $W_M$, and that $W_M$ is
an affine image of $W$.

In the same way one can compute the Rauzy fractal for the case $n=3$
in (\ref{substn}). The resulting pattern is shown in
Fig.~\ref{wkp9}, right. Note, that we did not claim that one of the
original substitution sets $V$ is a model set until now. In order to
prove this, one needs to prove $\mu(W)>0$, $\mu(\partial W)=0$ and
$\dens(V)=\dens(\{ \pi_1({\bs x}) \, | \, {\bs x} \in  \Lambda,
\pi_2({\bs x}) \in W \})$.  

\subsection{The dual tiling} \label{sec:dual}
\begin{figure}[h] 
\epsfig{file=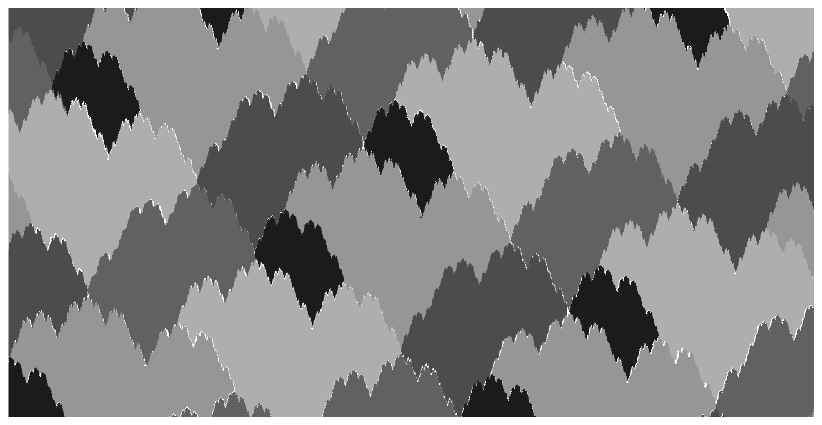} \epsfig{file=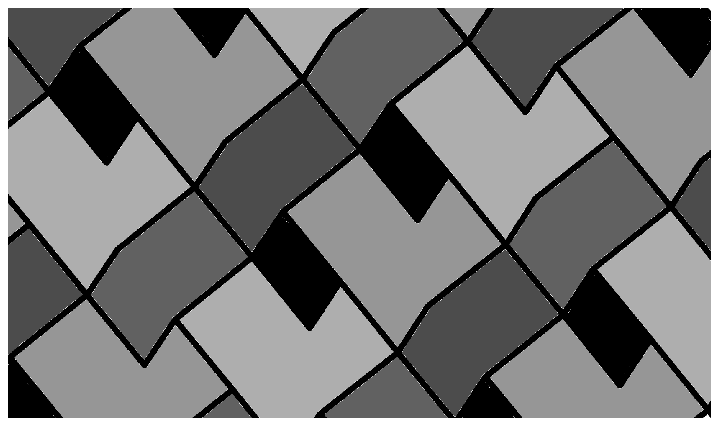}
\caption{ \label{fractilg}
Left: A patch of the tiling arising from (\ref{dualsubst}). Tiles of
type $W_S$ are black, the other tiles come in two slightly different
shadings, such that the boundary between tiles of the same type
remains visible. Right: A polygonalized version. Tiles and colours
correspond to the tiling on the left.}  
\end{figure}

Recall that $Q$ in (\ref{vstern}) is a contraction. Therefore $Q':=Q^{-1}$
is an expansion. Applying $Q'$ to (\ref{vstern}) and replacing 
$V_i^{\ast}$ with its closure $W_i$ ($i \in \{S,M,L\}$)  
yields the following equation system: 
\begin{equation} \label{dualexp}
\begin{array}{lcl}
Q' W_L & = & W_L \cup W_M \cup W_S \cup W_L - (\lambda_2,
\lambda_3)^T + (1,1)^T  \\
Q' W_M & = & W_L - (1,1)^T  \cup W_M - (1,1)^T
\cup W_S - (1,1)^T  \\ 
Q' W_S & = & W_M - (\lambda_2, \lambda_3)^T + (1,1)^T 
\end{array}
\end{equation}
This defines the substitution $\sigma'$ in $\E^2$ with three prototiles
$W_L,W_M,W_S$:
\begin{equation} \label{dualsubst}
\begin{array}{lcl}
\sigma' (\{W_L\}) & = & \{ W_L , W_M , W_S , W_L - (\lambda_2,
\lambda_3)^T + (1,1)^T \} \\
\sigma' (\{W_M\}) & = & \{ W_L - (1,1)^T  , W_M - (1,1)^T
, W_S - (1,1)^T \} \\ 
\sigma' (\{W_S\}) & = & \{ W_M - (\lambda_2, \lambda_3)^T + (1,1)^T \} 
\end{array}
\end{equation}
In contrary to Def.~\ref{substdef}, the expansion by $Q'$ is not a
similarity, but an affinity. Anyway, these substitution rules will
yield a tiling of $\E^2$. A part of it is shown in
Fig.~\ref{fractilg}. (To be precise, we need to show that there
will be no overlaps and no gaps, but this is not essential for the
following, since we will consider the corresponding substitution 
point set, or a modified version where all tiles are polygons.) 

Now, we will prove that the dual tiling is a model set.
In analogy to the construction of $V$ out of the original substitution
tiling arising from (\ref{subst2}), we choose in every prototile a 
{\em control  point}. A good choice is (cf.~(\ref{vstern})): 
$c_L:=0 \in W_L$, $c_M:= Qc_L- (\lambda_2,\lambda_3)^T =
(\lambda_2,\lambda_3)^T \in W_M$, $c_S:= Qc_M
-({\lambda_2}^2,{\lambda_3}^2)^T + (\lambda_2,\lambda_3)^T 
\in W_S$. It is not too hard to check that the set $\{W_L, W_M 
-({\lambda_2}^2,  {\lambda_3}^2)^T + 2 (\lambda_2,\lambda_3)^T + (1,1)^T\}$ 
gives rise to a point set $V'$ such that $\sigma'(V')=V'$, and $V'$ is
a Delone set. To be precise, let $V_0':=\{c_L, c_M-({\lambda_2}^2,
{\lambda_3}^2)^T + 2 (\lambda_2,\lambda_3)^T + (1,1)^T\}$, and
$\sigma'$ be defined for points as well as for tiles in the canonical way.
Then $V_0' \subset \sigma'(V_0')$, and inductively $(\sigma')^k (V_0') 
\subset  (\sigma')^{k+1} (V_0')$. Therefore, 
\begin{equation} \label{vstrich}
V' = \bigcup_{k \ge 1} (\sigma')^k (V_0') 
\end{equation} 
is well defined, and it holds $\sigma'(V')=V'$. 

{\em Remark:} As a first idea, one would try to start with the set 
$\{ c_L \}$ instead of $V_0'$. In this case, the property 
$(\sigma')^k(\{ c_L \}) \subset (\sigma')^{k+1}(\{ c_L \})$ is
fulfilled, but the resulting point set would be no Delone set, since
it possesses arbitrary large holes. 

By the construction, it holds $V \subset \langle (1,1)^T, (\lambda_2,
\lambda_3)^T , ({\lambda_2}^2,{\lambda_3}^2)^T \rangle_{\Z} $. In
particular, this shows a 
curious property of $V'$: If one coordinate of some ${\bs z}=(x,y)^T \in
V'$ is known, say, $x=a + b \lambda_2 + c {\lambda_2}^2$ ($a,b,c \in \Z)$,
then the other coordinate $y$ is known, too: $y = a  + b \lambda_3 + c
{\lambda_3}^2$. Another unusual property is that each prototile occurs
in just one orientation. This is not true for most known examples of
nonperiodic substitution tilings in dimension larger than one. 
    
In a similar way as above
(cf.~(\ref{vgleich}))  we obtain from (\ref{dualsubst}) the following
equations, where $V'=V'_S \cup V'_L \cup V'_L$ and $V'_i$ is the set
of all points of type $c_i$ in $V'$ ($i \in \{S,M,L\}$).    
\begin{equation} \label{dualgleich}
\begin{array}{lcl}
V'_L & = & Q' V'_L \cup Q' V'_M \cup Q' V'_L - (\lambda_2,
\lambda_3)^T + (1,1)^T  \\
V'_M & = & Q' V'_L - (\lambda_2,\lambda_3)^T  \cup Q' V'_M - (\lambda_2,
\lambda_3)^T \cup Q' V'_S   \\ 
V'_S & = & Q' V'_L + (\lambda_2, \lambda_3)^T - 2 ({\lambda_2}^2,
{\lambda_3}^2)^T \cup  Q' V'_M + (\lambda_2, \lambda_3)^T - 2
({\lambda_2}^2, {\lambda_3}^2)^T 
\end{array}
\end{equation}
In order to show that $V'$ is a model set, we have to determine the
lattice $\Lambda'$ and the window set $W'$ for $V'$. By the
construction of $V'$ follows $\Lambda'=\Lambda$. Since $0 \in V'$, we
know that every ${\bs x} \in V'$ is of the form
$a(1,1)^T+b(\lambda_2,\lambda_3)^T+c({\lambda_2}^2,{\lambda_3}^2)^T$. 
Thus our new star map is just the inverse of the original one in
(\ref{ast}), namely 
\[ \star: \pi_2(\Lambda) \to \E^1, \quad {\bs x}^{\star}=
\pi_1(\pi_2^{-1}({\bs x}))=a+b \lambda + c \lambda^2.  \] 
Applying $\star$ to (\ref{dualgleich}) yields:
\begin{equation} \label{dualglstar}
\begin{array}{lcl}
{V'_L}^{\star} & = & \lambda^{-1} {V'_L}^{\star} \cup \lambda^{-1}
{V'_M}^{\star} \cup \lambda^{-1} {V'_L}^{\star} - \lambda + 1  \\ 
{V'_M}^{\star} & = & \lambda^{-1} {V'_L}^{\star} - \lambda \cup \lambda^{-1}
{V'_M}^{\star} - \lambda \cup \lambda^{-1} {V'_S}^{\star} \\  
{V'_S}^{\star} & = & \lambda^{-1} {V'_L}^{\star} + \lambda - 2 \lambda^2 \cup
\lambda^{-1} {V'_M}^{\star} + \lambda - 2 \lambda^2  
\end{array}
\end{equation}
Since all occurring maps are contractions, Theorem \ref{hutch} applies,
and the unique compact solution is 
\[ \cl({V_L}^{\star}) = [-\lambda,0] = L, \; \cl({V'_M}^{\star}) =
[\lambda-\lambda^2,-\lambda] = M, \;  \cl({V'_S}^{\star}) = [-2
\lambda^2 + 1, -2 \lambda^2 + \lambda] = S. \]
The multiplication of (\ref{dualglstar}) by $\lambda$ yields a
distinct geometric realization of our original substitution
(\ref{subst2}). So, the dual tiling of the dual tiling is the one we
started with.  In particular, the window set of $V'$ is $W'=S \cup M
\cup L = [-2 \lambda^2 + 1, -2 \lambda^2 + \lambda] \cup
[\lambda-\lambda^2,0]$. Therefore $W'$ is compact, 
\begin{equation} \label{volw}
\mu(W')= \lambda^2-1 =  7.290859374\ldots >0
\end{equation}
and $\mu(\partial W')=0$. All we are
left with is to compute $\dens(V')$ and use Theorem \ref{densvol}.
If the tiles $W_i$ of the substitution $\sigma'$ in (\ref{dualsubst})
would be polygons, this would be an easy task, using the remark after
Theorem \ref{perron}. Therefore we replace the fractal tiles
$W_S,W_M,W_L$ by appropriate polygons $P_S,P_M,P_L$. Let  
\[ \begin{array}{rclrcl}
a & := & ({\lambda_2}^2 -1 ,  {\lambda_3}^2 -1 )^T , & 
b & := & ({\lambda_2}^2, {\lambda_3}^2 )^T,  \\ 
c & := &  ({\lambda_2}^2 - \lambda_2 ,  {\lambda_3}^2 -\lambda_3 )^T, &
d & := & ({\lambda_2}^2 - \lambda_2 -1,  {\lambda_3}^2 -\lambda_3 -1)^T, \\ 
e & := & (2{\lambda_2}^2 - \lambda_2 -1,  2{\lambda_3}^2 -\lambda_3 -1)^T, & 
f & := & (2{\lambda_2}^2 - 1,  2{\lambda_3}^2 -1)^T,  \\
g & := & (\lambda_2,\lambda_3)^T, & 
h & := & ({\lambda_2}^2 + \lambda_2 -1,  {\lambda_3}^2 +\lambda_3 -1)^T. 
\end{array} \]

\begin{figure}[h]
\epsfig{file=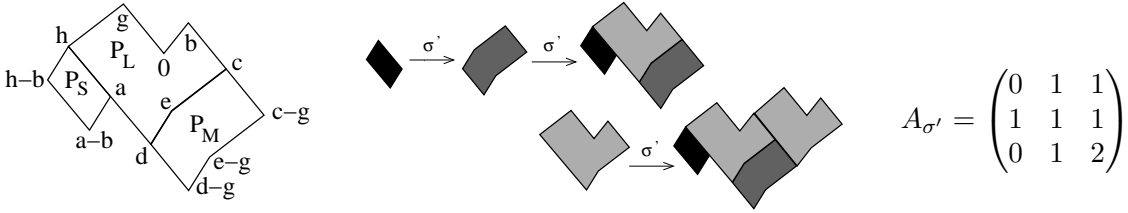} \hspace{5mm}
\begin{minipage}[b]{34mm}
$A_{\sigma'}= \begin{pmatrix} 0 & 1 & 1 \\
  1 & 1 & 1 \\ 0 & 1 & 2 \end{pmatrix}$ 
\vspace*{5mm}
\end{minipage}
\caption{Left: The construction of the polygonal prototiles. Middle:
A sketch of the substitution for them. Right: The corresponding
substitution matrix $A_{\sigma'}$. Note that $A_{\sigma'} =
(A_{\sigma})^T$. \label{polytiles}}
\end{figure}

Let $\conv({\bs x}_1,\ldots,{\bs x}_n)$ denote the convex hull of
${\bs x}_1,\ldots,{\bs x}_n$. Then (see Fig.~\ref{polytiles})
\begin{equation}
\begin{split}
P_S & := \conv(a-b,a,h,h-b), \\ 
P_M & := \conv(e-g,e,d,d-g) \, \cup \, \conv(c-g,c,e,e-g), \\ 
P_L & := \conv(0,e,c,b) \, \cup \, \conv(0,g,h,d,e)  
\end{split}
\end{equation}
It is quite simple (but lengthy) to show that by replacing $W_i$ with
$P_i$ in (\ref{dualsubst}), the substitution $\sigma'$ leads to a tiling,
without any gaps or overlaps.  
(Essentially, one uses the definition of a substitution tiling in
Definition \ref{substdef} and the primitivity of the substitution
matrix. To show that there are no gaps and no
overlaps it suffices to consider all vertex configurations in the
first $n$ substitutions $\sigma'(T),(\sigma')^2(T), \ldots,
(\sigma')^n(T)$ for appropriate $n$ and some prototile $T$. So, in any
certain case this is a finite   
problem. Here $n=4$ and $T=P_L$ will do.) It is also quite simple (but
lengthy again) to compute the areas of the prototiles
$P_S,P_M,P_L$. Unfortunately, we cannot use the fact that the left
eigenvector of the corresponding substitution matrix reflects the
ratio of the areas, since the substitution for the polygonal tiles
does not preserve the shape any longer. E.g., $\sigma'(W_S)$ is congruent
to $Q'W_S$, but $\sigma'(P_S)$ is {\em not} congruent to $Q'P_S$. 
Altogether, we obtain  
\begin{equation}
\begin{split} 
  \mu(P_S) & = |\lambda_2| \lambda_3 (\lambda_3-\lambda_2) = \frac{s_1}{s_4}
  \left(  \frac{s_2}{s_4} + \frac{s_1}{s_2} \right), \\ 
  \mu(P_M) & = \lambda_3-\lambda_2 = \frac{s_2}{s_4} + \frac{s_1}{s_2}, \\ 
  \mu(P_L) & = {\lambda_2}^2 - {\lambda_3}^2 + 2 |\lambda_2| + 2 \lambda_3 -
  |\lambda_2| {\lambda_3}^2 - \lambda_3 {\lambda_2}^2 = 1 + \frac{s_4
  s_1}{{s_2}^2}, \\
\end{split}
\end{equation}
where the last equality requires an excessive use of (\ref{wformel}).
The vector $(\mu(P_S),\mu(P_M),\mu(P_L))$ indeed is a left
eigenvector of the substitution matrix $A$ (cf.~Fig.~\ref{polytiles})
corresponding to the PF eigenvalue $\lambda$. A right eigenvector
corresponding to $\lambda$ is $(s_2,s_3,s_4)^T$. Therefore,  
the normalized right eigenvector of $A$ --- which entries, by Theorem
\ref{freqvol}, are the
relative frequencies $p_i$ of tile types $P_i$ ($i \in \{S,M,L\}$) in
the tiling --- is $(p_S,p_M,p_L)^T := \frac{1}{s_2 + s_3 + s_4}
  (s_2,s_3,s_4)$. Since every point of $V'$ corresponds to exactly one
tile in the polygonal tiling, by Prop. \ref{densfreqvol} the relative
density of $V'$ is 
\begin{equation} \label{densv}
\begin{split}
\dens(V') & = 
\left( (\mu(P_S),\mu(P_M),\mu(P_L))\cdot(p_S,p_M,p_L)^T \right)^{-1} \\
& = (s_2+s_3+s_4) \left( \frac{s_1{s_2}^2}{s_4^2} + \frac{s_1^2}{s_4}
+ \frac{s_2 s_3}{s_4}+\frac{s_1 s_3}{s_2}+ s_4
+ \frac{s_1 {s_4}^2}{{s_2}^2} \right)^{-1} \\
& = 0.8100954858 \ldots
\end{split}
\end{equation}
In order to apply Theorem \ref{densvol} we need to compute $\det
\Lambda$ (cf.~(\ref{mlambda})). A rather lengthy computation (or the 
usage of a computer algebra system, here: Maple) shows 
$\det \Lambda = \frac{27}{64 {s_1}^2 {s_2}^2  {s_4}^2}$. Using
elementary trigonometric formulas yields
\[ 
\begin{split}
s_1 s_2 s_4 & = \frac{1}{4} (s_3+s_4 - s_2 - s_1 ) \\
& = \frac{1}{4} \left( s_3 + 2 \cos \left( \frac{3 \pi}{9} \right)
  s_1 -s_1 \right) = \frac{1}{4} s_3 = \frac{\sqrt{3}}{8}, 
\end{split}
\] 
and therefore $\det \Lambda = 9$. Putting together this with
(\ref{volw}) and (\ref{densv}) in view of Theorem \ref{densvol} we
obtain \[ \mu(W')/\dens(V') = 9 =\det(\Lambda). \]
The leftmost equality requires again a very lengthy computation, or
the usage of Maple. This completes the proof, and we obtain the
following theorem. 
\begin{thm} \label{dersatz}
The set $V'$ in (\ref{vstrich}), arising from the substitution $\sigma'$
in (\ref{dualsubst}), is a model set. 
\end{thm}
One may ask why we do not state that the original set $V$ arising from
the substitution (\ref{subst2}) is a model set. The reason is that
there is one gap, namely, the value of $\mu(W)$. For example, one needs
to show that the tiling by the fractal tiles $W_S,W_M,W_L$ is not
overlapping. This seems to be clear from the pictures, but we did not
give a rigorous proof here. Indeed, this will require some more work,
and we refer to further publications. Anyway, there are good reasons
to assume that the areas of the 
fractal tiles $W_S,W_M,W_L$ are equal to the areas of the polygonal
tiles $P_S,P_M,P_L$, resp. Proving this will complete the
proof that the original set $V$ is a model set, too. 

\section{Further Remarks} \label{sec:remarks}
The question whether a given structure is a model set is motivated for 
example by the question for the diffraction property of this structure. 
It is known that model sets have a pure point diffraction image,
without any continuous part. For details, cf.~\cite{bmrs}.

We did not emphasize it throughout the text, but all tilings considered
here are nonperiodic. For the one--dimensional tilings this is a
consequence of the fact that the relative frequencies of the tiles are
irrational numbers. If a one--dimensional tiling is periodic, then the
relative frequencies are rational numbers. Anyway, all occurring
two--dimensional tilings are nonperiodic, too. This can be shown with
standard tools from the theory of nonperiodic substitution tilings
(cf.~statement 10.1.1 in \cite{gs}). 

An important point in proving Theorem \ref{dersatz} is that we know
the occurring values in terms of trigonometric expressions, thus
it is possible to make exact computations of $\mu(W')$, $\dens(V')$,
$\det(\Lambda)$ etc. In general, the fractal nature of the occurring 
window sets may forbid this. But if all occurring sets are polytopes
the things become easier. In the following, some examples are listed
where the duality principle provides short and simple proofs that
these considered substitution tilings are model sets.

{\bf Fibonacci sequences:} The symbolic substitution 
$S \to L,\, L \to LS$ gives rise to the widely examined Fibonacci
sequences. (Not to be confused with the Fibonacci numbers
$1$,$1$,$2$,$3$,$5$,$8$,$13$,$\ldots$, sometimes also called Fibonacci
sequence. If one considers the first substitutions of $S$, then the
relation becomes clear: $L$, $LS$, $LLS$, $LSLLS$, $LSLLSLSL$, $\ldots$). 

A geometric representation as a substitution tiling
(Fibonacci tiling) is for
example given by $L:=[0,1],S=[1, \tau]$ ($\tau$ the golden ratio)
and \[ \sigma (\{S\}) = \{ L+\tau \}, \; \sigma( \{L\} ) = \{L,S\}. \]
By the standard construction (cf.~Section \ref{sec:modset}), we obtain
for some $V$ which fulfills $\sigma(V)=V$ (namely, $V= \bigcup_{k \ge
  0} \sigma^{2k}(\{ L, L-1\})$, up to one point this will do) the
equations 
\[ \begin{split}
V_L^{\ast} & = -\tau^{-1} V_L^{\ast} \cup -\tau^{-1} V_S^{\ast} \\
V_S^{\ast} & = -\tau^{-1} V_L^{\ast} + 1
\end{split}, \]
to be read as equations in $\E^1$. The unique compact solution is 
\[ \cl(V_L^{\ast}) = W_L = [-1, \tau^{-1}] , \;    
\cl(V_S^{\ast}) = W_S = [\tau^{-1}, \tau]. \]
The dual tiling is given by the following substitution:
\[ \begin{split} 
\sigma'(\{W_L\}) & = \{W_L,W_S\} \\
\sigma'(\{W_S\}) & = \{W_L - \tau \} 
\end{split} \]
Note that the substitution factor is $-\tau$. Nevertheless, the
resulting family of substitution tilings for $\sigma'$ is equal to
these for $\sigma$. That means the dual tilings of Fibonacci tilings
are Fibonacci tilings again; or shortly: Fibonacci tilings are
self--dual. Moreover, a short computation yields that the equation in
Theorem \ref{densvol} is fulfilled, and therefore the Fibonacci
sequences are model sets. 

\begin{figure}[h]   
\epsfig{file=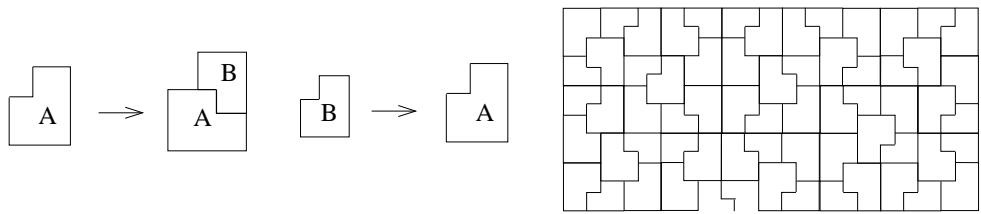}
\caption{\label{ammch} Left: The substitution for the Ammann
  chair. Right: A part of a corresponding substitution tiling.}
\end{figure}

{\bf Golden triangles and Ammann chair:}  
The substitution for the golden triangles (cf.~Fig.~\ref{inflbsp})
defines a model set, if one replaces every triangle with the vertex
at its right angle. This was essentially shown in \cite{dvo} in
detail\footnote{The authors used {\em all} vertices in the tiling
  instead of one vertex of each triangle, but their result is easily
  transferred to the latter case.}. It turns out that the golden
triangle tilings are dual to the so--called Ammann chair tilings (or
'A2 tilings', cf.~\cite{gs}). A sketch of the substitution of the
latter is shown in Fig.~\ref{ammch}, together with a part of the
tiling. The substitution factor and the substitution matrix are equal
to those of the golden triangle substitution. All occurring tiles, and
all occurring window sets are polygons. Therefore, it is simple to
compute the areas and the densities, thus proving --- with the help of
Theorem \ref{densvol} --- that both substitution tilings (resp.~the
corresponding point sets) are model sets. 

The duality of the Ammann chair tilings and the golden triangle
tilings was already indicated in \cite{gel} briefly. But a proof of
this needs some care. For example, one has to pay attention to the
fact, that each prototile occurs in four different orientations, not
just in one, as it is true for the tilings in Fig.~\ref{fractilg}.  

\begin{figure}[h]   
\epsfig{file=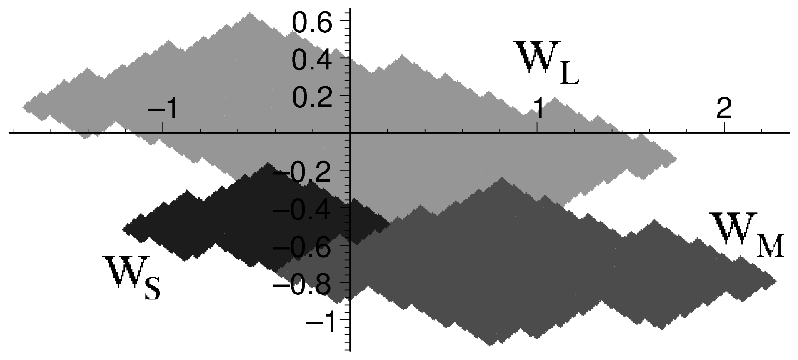} \hfill \epsfig{file=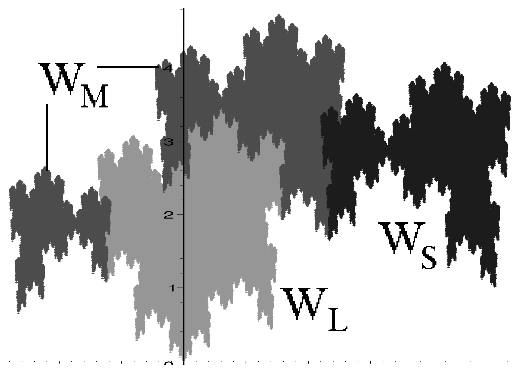}
\caption{\label{mehrrosi} 
Two more Rauzy fractals, the left one arising from the substitution
(\ref{substperm}), the right one arising from the substitution
(\ref{substtransp}) }
\end{figure}

{\bf Variants of the tilings in Section \ref{sec:zwo}:} 
The tilings considered in the last section are a rich source of other 
substitutions, or other Rauzy fractals. If we restrict to the case
$n=4$, then the underlying lattice $\Lambda$ will always be the same
(namely, as in (\ref{mlambda})). One possible modification is to
permute the letters in the substitution. E.g., the substitution
\begin{equation} \label{substperm}
 S \to ML, \; M \to MLS, \; L \to LML 
\end{equation}
instead of (\ref{subst2}) yields the Rauzy fractal shown in
Fig.~\ref{mehrrosi} (left). Another variant is to use the
matrix $A_{\sigma'}$ of Fig.~\ref{polytiles} (the transpose of the
original substitution matrix $A_{\sigma}$) as a substitution matrix
for a one--dimensional substitution with substitution factor
$\lambda$. One possibility is the substitution  
\begin{equation} \label{substtransp} 
S \to M, \; M \to SML, \; L \to SMLL, 
\end{equation}
where $S,M,L$ are intervals of length $\frac{s_1 s_2}{s_4}, s_2,s_4$,
resp. This substitution yields the Rauzy fractal shown in
Fig.~\ref{mehrrosi} (right). Note that in this case the set $W_M$ is not
connected, it consists of two parts which are similar to each other. 
The substitution tiling arising from this Rauzy fractal (resp.~from
the corresponding equation system, multiplied by $Q'$,
cf.~(\ref{dualexp})) would possess a prototile, namely, $W_M$, which is
disconnected. Note, that we did not exclude disconnected
(proto--)tiles in our definition of (substitution) tilings.

Combining these two modifications leads to even more examples, all of
them based on the lattice $\Lambda$. All the two--dimensional tilings
are model sets, which is easily proven by using the results of Section
\ref{sec:zwo}.   
From the case $n=3$ one can also construct different tilings by
permuting the letters. In contrary, all variants of the case $n=2$
will yield Fibonacci sequences. Anyway, in the cases $n=3, n=4$ it may
be interesting to examine the different possibilities with respect to
common properties of the corresponding Rauzy fractals.

\end{document}